\documentclass[11pt]{article}
\usepackage[T1]{fontenc}
\usepackage[margin=1.03in]{geometry}
\usepackage{amsmath,amsthm,mathtools}
\usepackage{newtxtext,newtxmath}
\usepackage{microtype}
\usepackage[hidelinks]{hyperref}
\hypersetup{pdftitle={Adversarial Wronskians: A proof of the Kung--Traub conjecture},pdfauthor={Matthew J. Colbrook},pdfsubject={Kung--Traub conjecture, local order, and polynomial Wronskians}}
\usepackage[nameinlink,noabbrev,capitalize]{cleveref}
\usepackage{enumitem}
\usepackage{tikz}
\setlist{nosep}
\numberwithin{equation}{section}
\newtheorem{theorem}{Theorem}[section]
\newtheorem{lemma}[theorem]{Lemma}
\newtheorem{corollary}[theorem]{Corollary}
\theoremstyle{definition}

\theoremstyle{remark}
\newtheorem{remark}[theorem]{Remark}
\newcommand{\R}{\mathbb R}
\newcommand{\C}{\mathbb C}
\newcommand{\Pn}{\mathcal P}
\renewcommand{\Wr}{\operatorname{Wr}}
\newcommand{\dist}{\operatorname{dist}}

\newcommand{\diag}{\operatorname{diag}}

\newcommand{\norm}[1]{\left\lVert#1\right\rVert}
\newcommand{\abs}[1]{\left\lvert#1\right\rvert}
\newcommand{\ball}{\overline B}

\title{Adversarial Wronskians: A proof of the Kung--Traub conjecture}
\author{Matthew J. Colbrook\thanks{Department of Applied Mathematics and Theoretical Physics, University of Cambridge, Wilberforce Road, Cambridge CB3 0WA, United Kingdom. Email: \href{mailto:m.colbrook@damtp.cam.ac.uk}{\texttt{m.colbrook@damtp.cam.ac.uk}}.}}
\date{19 September 2026}
\begin{document}
\maketitle
\vspace{-4mm}

\begin{abstract}
We prove the Kung--Traub conjecture: a deterministic stationary root-finding method without memory, using at most $n$ real scalar function or derivative evaluations per iteration, cannot have local order greater than $2^{n-1}$ on every real analytic function with a simple zero. For $n\geq2$, the bound is attained by the methods of Kung and Traub. The local estimate is required at every sufficiently close starting point, with constants that may depend on the function. Evaluation points and derivative orders may be chosen adaptively, and the rules for choosing them and processing the data are unrestricted. The proof uses a localisation property of polynomial Wronskians to obtain a lower bound at each fixed scale. A diagonal construction then gives one entire function, with derivative bounded above and below on the real line, for which every larger order fails. We also obtain a sharper bound when evaluations are made in prescribed groups, and prove the corresponding result for complex scalar evaluations. The results of this paper are also formally verified in Lean.
\end{abstract}

\section{Introduction}

We prove the conjecture of Kung and Traub \cite[Conjecture~8.1]{KT}, proposed in 1974, that a root-finding iteration without memory using $n$ function or derivative evaluations has local order at most $2^{n-1}$. Newton's method attains this bound with two evaluations. Kung and Traub constructed multipoint methods attaining it for every $n\geq2$; the result below establishes the matching upper bound.

To state our main result, we first describe the precise setup. An input is a real analytic function $f$ on an open interval, that has a specified simple zero $\alpha_f$. The root $\alpha_f$ specifies where local order is measured; its value is not supplied to the algorithm. One execution of an algorithm starts with a real number $x$. At step $j$ it chooses a real location $z_j$ and an integer $\nu_j\geq0$ from $x$ and the previous answers, and receives
\[
 Y_j=f^{(\nu_j)}(z_j).
\]
After at most $n$ such observations it returns a real number $T_f(x)$ depending only on $x$ and these answers. No continuity or measurability of the rules choosing $z_j,\nu_j$, or the final output is imposed. The rules are fixed and are the same for every function and every iteration. In particular, there is no information retained from earlier iterations other than the current starting point. An execution is defined only when its query points lie in that interval and its stage rules are defined.\footnote{An algorithm that stops early can be padded with zero linear observations, which have no associated query location. Locally defined stage rules are extended arbitrarily on undefined transcripts. If the original method is defined for all starting points sufficiently close to a simple zero, the extension agrees with it there.} The counterexamples to local order constructed below are entire, so their domains impose no restriction on the query points.

We say that the method has universal local order at least $p$ if, for every such input, it is defined for every sufficiently close starting point $x\ne\alpha_f$ and satisfies
\begin{equation}\label{eq:localorder}
 \abs{T_f(x)-\alpha_f}=O_f\bigl(\abs{x-\alpha_f}^{p}\bigr)
 \quad\text{as }x\rightarrow\alpha_f.
\end{equation}
The constants and neighbourhood may depend on $f$; the estimate is required for every $x\ne\alpha_f$ in that neighbourhood.


\begin{theorem}[Local order bound]\label{thm:main}
Let $n$ be a positive integer, and let $T$ be an algorithm of the kind just described. For every real $p>2^{n-1}$, there exist a function $f$ entire on $\C$, real-valued on $\R$, and a sequence $(x_s)_{s=1}^{\infty}$ of real starting points such that
\[
 \frac34\leq f'(t)\leq\frac54\quad(t\in\R),\qquad
 x_s\rightarrow\alpha\quad(s\to\infty),
\]
where $\alpha$ is the unique real zero of $f$, $x_s\ne\alpha$ for every $s\geq1$, and
\begin{equation}\label{eq:failure}
 \frac{\abs{T_f(x_s)-\alpha}}{\abs{x_s-\alpha}^{p}}
 \rightarrow\infty\quad(s\to\infty).
\end{equation}
Hence, no such method can satisfy \eqref{eq:localorder} with $p>2^{n-1}$ for every real analytic function with a simple zero.
\end{theorem}

The definition of order in \cite[Section~2]{KT} implies \eqref{eq:localorder}, and every iteration considered there belongs to the present information model. Thus \cref{thm:main} proves the Kung--Traub conjecture. Together with the matching construction in \cref{sec:sharpness}, it establishes the optimal order $2^{n-1}$ for $n\geq2$. The theorem uses real observations; \cref{sec:complex} proves the corresponding result for complex observations on the unrestricted class of holomorphic functions.

\begin{remark}
Our argument permits arbitrary deterministic processing, evaluations at unrestricted real points, and derivatives of any nonnegative integer order. In particular, the derivative data need not consist of consecutive derivatives at a point. The counterexample can be chosen entire and strictly increasing, with derivative uniformly close to one. The obstruction therefore persists for functions with a unique, well-conditioned real zero.
\end{remark}

A more precise bound holds when evaluations are made in prescribed groups.

\begin{corollary}[Grouped evaluations]\label{cor:batch}
Let $T$ be an algorithm as above. Suppose each execution has $k\geq1$ groups of observations, of fixed positive integer sizes $\ell_1,\dots,\ell_k$. All queries in group $j$ may depend on the earlier groups' answers, but not on answers within group $j$. Then its universal local order is at most
\begin{equation}\label{eq:batchbound}
 \ell_1\prod_{j=2}^{k}(\ell_j+1).
\end{equation}
The derivative orders within each group are unrestricted.
\end{corollary}

For a total of $n=\sum_{j=1}^{k}\ell_j$ observations, a single group gives the upper bound $n$, whereas $n$ groups of size one recover $2^{n-1}$. If every group has size $\ell$, the exponent is $\ell(\ell+1)^{k-1}$. This exponent coincides with the upper bound obtained by Casulli and Trigiante for their parallel multipoint procedures without memory \cite{CT}.

Lean verification for the results of this paper can be found here:
\begin{center}
\url{https://github.com/MColbrook/Kung-Traub_conjecture}.
\end{center}
The formal verification is also registered on Palomar:
\begin{center}
\url{https://palomar-registry.org/entry.html?id=PALOMAR-2026-09-20-000006&version=1}.
\end{center}

\subsection{Earlier results}

The conjecture extends Traub's earlier conjecture for two and three evaluations \cite{Traub1972}; see also \cite[Section~8]{KT}. Kung and Traub proved the bound $2^{n-1}$ for Hermite interpolatory iterations and a bound $2^n$ for iterations using function values alone \cite[Section~7]{KT}. They subsequently established the conjectured bound for rational iterations using two evaluations \cite{KTtwo}.

Wo\'zniakowski introduced the order of information for stationary operator iterations under conditions on the positions of the iterates \cite{Woz1975}. In this framework, he proved the bound for Hermitian information \cite{Woz}. Here Hermitian information means that the data at each point comprise the function value and all derivatives up to a prescribed order. The product in \eqref{eq:batchbound} also occurs in that setting. Wasilkowski obtained further cases through conditions on the associated Birkhoff interpolation problems \cite{Wasil}. The proof below does not require a prescribed derivative interpolation problem to be uniquely solvable.

The construction of one function preserving earlier observations has an antecedent in Brent, Winograd, and Wolfe \cite{BWW}. Their results concern convergence along generated sequences, with a different allowance for derivative information and memory. Casulli and Trigiante obtained order bounds for classes of multipoint procedures with and without memory \cite{CT1977a}. For methods with memory, Werschulz obtained an upper bound $2^n$ for multipoint Hermitian information with $n$ new evaluations per iteration \cite{Werschulz}.

\subsection{Outline of the proof}

The proof constructs functions that give exactly the same answers to the algorithm but have different roots. A deterministic algorithm must return the same value on such inputs, so its error on at least one of them is at least half the distance between their roots. We first prove a quantitative separation at each scale, and then obtain a single entire function for which the local estimate fails along a sequence of starting points. The mechanism is illustrated in \cref{fig:proof-mechanism}.

\begin{figure}[!t]
\centering
\begin{tikzpicture}[x=1cm,y=1cm,font=\normalsize,>=stealth]
\definecolor{proofblue}{RGB}{35,72,103}
\path[use as bounding box] (0,-.98) rectangle (15.5,7.7);

\node[anchor=west,font=\normalsize\bfseries] at (0,7.35)
  {(a) Choose the next answer};
\node at (4.65,6.8) {Attainable roots before the answer};
\fill[black!8] (3,6.17) circle (.46);
\fill[black!8] (3,5.53) circle (.46);
\draw[black!40] (3,6.17) circle (.46);
\draw[black!40] (3,5.53) circle (.46);
\fill[black!8] (6.2,5.85) circle (.53);
\draw[black!40] (6.2,5.85) circle (.53);
\draw[->,black!50] (.35,5.85) -- (15.1,5.85) node[right,text=black] {$\R$};
\draw[line width=.85pt] (.75,5.85) -- (14.55,5.85);
\draw (.75,5.73) -- (.75,5.97);
\draw (14.55,5.73) -- (14.55,5.97);
\draw[line width=.8pt] (2.9,6.07) -- (3.1,6.27);
\draw[line width=.8pt] (2.9,6.27) -- (3.1,6.07);
\draw[line width=.8pt] (2.9,5.43) -- (3.1,5.63);
\draw[line width=.8pt] (2.9,5.63) -- (3.1,5.43);
\fill (6.2,5.85) circle (2pt);
\draw[proofblue,line width=4pt] (9.7,5.85) -- (12.7,5.85);
\fill[proofblue] (11.2,5.85) circle (2.8pt);
\draw[proofblue,thin] (11.2,6.04) -- (11.2,6.34);
\node[proofblue,above] at (11.2,6.34) {Retained root interval};
\node[proofblue,below] at (11.2,5.64) {$\alpha_*$: chosen root};
\node[align=center,font=\small] at (3,4.77)
  {Wronskian zeros};
\node[align=center,font=\small] at (6.2,4.77)
  {Next query location};
\node[anchor=west,font=\small,text=black!70] at (.75,4.32)
  {Grey discs: neighbourhoods to avoid};

\draw[black!25] (.1,3.93) -- (15.4,3.93);
\node[anchor=west,font=\normalsize\bfseries] at (0,3.58)
  {(b) After $n$ observations: identical data, separated roots};
\node at (3.25,2.87) {$T_{f_-}(x)=T_{f_+}(x)=y$};
\node[font=\small] at (3.25,2.36) {The algorithm returns the same value};
\draw[->,black!50] (7,2.5) -- (15.1,2.5);
\fill[proofblue] (8,2.5) circle (2.6pt);
\fill[proofblue] (11.7,2.5) circle (2.6pt);
\fill (14.2,2.5) circle (2.6pt);
\node[below,proofblue] at (8,2.35) {$\alpha_-$};
\node[below,proofblue] at (11.7,2.35) {$\alpha_+$};
\node[below] at (14.2,2.35) {$y$};
\draw[<->] (8,2.97) -- (11.7,2.97);
\node[above,font=\small] at (9.85,2.97) {Root separation};
\node at (7.75,1.52)
  {$\displaystyle\max\{|y-\alpha_-|,|y-\alpha_+|\}\geq\frac12|\alpha_+-\alpha_-|$};

\draw[black!25] (.1,1.02) -- (15.4,1.02);
\node[anchor=west,font=\normalsize\bfseries] at (0,.65)
  {(c) Preserve earlier answers and pass to one entire function};
\node at (4.05,-.03)
  {$f_0\longrightarrow f_1\longrightarrow\cdots\longrightarrow f_s\longrightarrow\cdots\longrightarrow f$};
\node[font=\small] at (4.05,-.62) {All earlier query answers remain unchanged};
\node at (11.6,-.03)
  {$\displaystyle\frac{|T_f(x_s)-\alpha|}{|x_s-\alpha|^p}\longrightarrow\infty$};
\node[font=\small] at (11.6,-.72)
  {$x_s\longrightarrow\alpha,\qquad p>2^{n-1}$};
\end{tikzpicture}
\caption{The proof in three steps. (a) Once the next query is fixed, its location and the relevant Wronskian zeros determine the grey exclusion discs. Choose an input with root $\alpha_*$ outside these discs and supply its answer. Restricting the parameter ball retains an interval of roots separated from the forbidden points. (b) After $n$ observations, two inputs still have identical data and roots $\alpha_-,\alpha_+$ separated by at least a positive constant times $\varepsilon^{2^{n-1}}$. Their common output $y$ has error at least half this separation on one input. (c) Repeat at decreasing scales while preserving all earlier answers. The entire limit violates the order-$p$ estimate for any fixed $p>2^{n-1}$.}
\label{fig:proof-mechanism}
\end{figure}

Fix a starting point $x$ at distance $\varepsilon$ from the root of a reference function. We consider perturbations $w(t)q(t-x)$, where $w$ is fixed and nonzero near the root, and $q$ has degree at most $n$. The multiplier and coefficient bounds keep the first derivatives of the perturbed functions bounded above and away from zero. Scaling the constant coefficient of $q$ by $\varepsilon$ then keeps the perturbed roots within a constant times $\varepsilon$ of $x$. Each scalar observation imposes at most one linear condition on the $n+1$ coefficients. Nonzero perturbations therefore remain after $n$ observations. This dimension count alone is insufficient: a perturbation that preserves the data might also be very small at the root. The essential step is to control how much the root can still vary.

Polynomial Wronskians provide this control. Our localisation lemma states that if a nonzero member of an $m$-dimensional polynomial space has $m$ zeros in a small disc, then the Wronskian of that space has a zero in a disc with the same centre and a uniformly comparable radius. Once the next query is chosen, the kernels of all observations up to that query are known before its answer is supplied. At each observation, the preceding step supplies an interval of attainable roots whose length is bounded below by a constant times a power of $\varepsilon$. The number of zeros of the nonconstant Wronskians is uniformly bounded. We can therefore choose the answer from a function whose root stays a fixed fraction of this interval's length away from all of them. A smaller ball of parameters preserves this separation. We consider every degree up to $n$, since the minimum degree of a polynomial in the final kernel is not known in advance.

This also explains the exponent. After $j$ observations, take a nonzero polynomial $q$ of minimum degree in the common kernel. A constant polynomial gives the required variation directly. Otherwise, the localisation lemma at successive rank increases bounds the ordered distances to the roots of $q$. The product formula contributes at most one factor from each observation, with the exponent of the root interval available before it. The scaled coefficient estimate and derivative bound then give a lower bound for root variation with the sum of all preceding exponents. Starting from one, the exponent remains one after the first observation and then doubles, giving $2^{n-1}$ after $n$ observations.

Finally, the constant in \eqref{eq:localorder} may depend on the function, so separate examples at each scale do not suffice. We construct a sequence of entire corrections, each vanishing to sufficiently high order at every earlier query point. All earlier answers, and hence all earlier outputs, are then preserved exactly, even for discontinuous rules. The corrections are summable on every complex disc and in the real $C^1$ norm. At each stage the next scale is chosen only after its lower-bound constant is fixed. The entire limit retains a simple, well-conditioned root and violates \eqref{eq:localorder} along the resulting sequence. The two parts of the construction are given in \cref{sec:finite,sec:entire}. The grouped bound in \cref{cor:batch} is proved in \cref{sec:groups}. The matching construction and the complex extension appear in \cref{sec:sharpness,sec:complex}, respectively.

\section{Polynomial preliminaries}\label{sec:polynomials}

\subsection{Polynomial Wronskians}\label{sec:wronskians}

Throughout, $\ball(a,r)$ denotes the closed Euclidean ball with centre $a$ and radius $r$, in the space indicated by the context. Let $\Pn_d$ be the space of complex polynomials in one variable of degree at most $d$. For an $m$-dimensional subspace $V\subseteq\Pn_d$ and a basis $p_1,\dots,p_m$, define
\[
 W_V(t)=\Wr(p_1,\dots,p_m)(t)
       =\det\bigl[p_i^{(k-1)}(t)\bigr]_{k,i=1}^{m}.
\]
Changing the basis multiplies $W_V$ by a nonzero constant, so its zero set and degree depend only on $V$. For background on polynomial Wronskians and the Wronski map, see \cite{EG}.

Choose a basis with distinct degrees
$0\leq d_1<\cdots<d_m\leq d$. The coefficient of the top power in the determinant is the product of the leading coefficients times
$\prod_{i<j}(d_j-d_i)$. Thus
\begin{equation}\label{eq:wrdegree}
 W_V\not\equiv0,\qquad
 \deg W_V=\sum_{i=1}^{m}d_i-\frac{m(m-1)}2
 \leq m(d-m+1).
\end{equation}
In particular, $W_V$ is a nonzero constant if and only if
$V=\Pn_{m-1}$.

The vanishing of $W_V$ at a point is equivalent to the existence of a nonzero member of $V$ with a zero there of multiplicity at least $m$; see also \cite[Theorem~A]{Dyakonov} for holomorphic functions. We need a uniform version in which the zeros may lie at distinct nearby points.

\begin{lemma}[Wronskian localisation]\label{lem:wr}
For $1\leq m\leq d$ there is a finite positive constant $C(d,m)$ such that the following holds. If an $m$-dimensional subspace $V\subseteq\Pn_d$ contains a nonzero polynomial with at least $m$ zeros in the closed disc $\ball(a,r)$, counted with multiplicity, then $W_V$ has a zero in
$\ball(a,C(d,m)r)$.
\end{lemma}
\begin{proof}
For $r>0$, an affine change of variable reduces the problem to $a=0,r=1$. Suppose there is no such constant. For each positive integer $\ell$, choose $V_\ell$ and a nonzero $h_\ell\in V_\ell$ such that $h_\ell$ has $m$ zeros in $\ball(0,1)$ but $W_{V_\ell}$ has no zero in $\ball(0,\ell)$.

Identify polynomials with their coefficient vectors in $\C^{d+1}$. Normalise $h_\ell$ to Euclidean norm one and extend it to an orthonormal basis of $V_\ell$. Compactness of the finite-dimensional Stiefel manifold gives a subsequence of these bases converging to an orthonormal basis of an $m$-dimensional subspace $V$. Its first member $h$ is nonzero.

The limit $h$ still has at least $m$ zeros in $\ball(0,1)$. Indeed, choose $m$ such roots of each $h_\ell$ with multiplicity and pass to convergent subsequences. Their monic degree-$m$ factors converge. Polynomial division by these monic factors shows that the quotient coefficients are bounded, so a further subsequence converges and the limiting factor divides $h$.

The Wronskians of the chosen bases converge coefficientwise to $W_V$, which is nonzero by \eqref{eq:wrdegree}. Any finite zero of $W_V$ would, by Rouch\'e's theorem on a sufficiently small circle, force nearby zeros of $W_{V_\ell}$ for large $\ell$. Thus $W_V$ has no finite zero and is a nonzero constant. Formula \eqref{eq:wrdegree} then forces $V=\Pn_{m-1}$, contradicting the $m$ zeros of its nonzero member $h$.

For $r=0$, choose the polynomial with an $m$-fold zero at $a$ as the first basis member. Its first $m$ derivatives starting at order zero vanish at $a$, so the first column of the Wronskian matrix vanishes there.
\end{proof}

We will need a uniform bound on the number of Wronskian zeros. Define
\begin{equation}\label{eq:Mn}
 M_n=\left\lfloor\frac{(n+1)^2}{4}\right\rfloor,
 \qquad H_n=nM_n+1.
\end{equation}
By \eqref{eq:wrdegree}, every nonconstant Wronskian of a subspace of $\Pn_d$, $d\leq n$, has at most $M_n$ distinct complex zeros. A zero-dimensional subspace is simply omitted, and the Wronskian of the full space $\Pn_d$ is constant. At each observation, we shall use one subspace for each degree $1,\dots,n$, together with one possible query location; $H_n$ bounds the resulting number of points.

\subsection{Coefficient estimates}\label{sec:coefficients}

For a polynomial $q(s)=\sum_{k=0}^{d}b_ks^k$ of degree $d\leq n$ in a centred variable, let $b=(b_0,\dots,b_d,0,\dots,0)\in\C^{n+1}$ be its coefficient vector, with zeros appended. We use the Euclidean norm of this vector.

The first estimate converts distances to the zeros into a lower bound for the polynomial's value. The second controls the change of normalisation when the constant coefficient is scaled by $\varepsilon$.

\begin{lemma}[Product and scaled-norm estimates]\label{lem:norm}
Let $q$ have degree $d\leq n$ and $\norm{b}_2=1$. If $\beta_1,\dots,\beta_d$ are the roots of $q$, counted with multiplicity, then for $|a|\leq1/2$,
\begin{equation}\label{eq:product}
 |q(a)|\geq 4^{-d}\prod_{k=1}^{d}\min\{1,|a-\beta_k|\}.
\end{equation}
For $\varepsilon>0$, let $D_\varepsilon$ be the $(n+1)\times(n+1)$ diagonal matrix $\diag(\varepsilon,1,\dots,1)$. For every fixed $R>0$, if $|a|\leq R\varepsilon\leq1/2$, then
\begin{equation}\label{eq:scaled}
 \frac{|q(a)|}{\norm{D_\varepsilon^{-1}b}_2}
 \geq \frac{1}{2R+2}\min\{\varepsilon,|q(a)|\}.
\end{equation}
\end{lemma}
\begin{proof}
The factorisation $q(s)=b_d\prod_k(s-\beta_k)$ gives
\[
 1\leq\norm{b}_1
 \leq |b_d|\prod_k(1+|\beta_k|)
 \leq 2^d|b_d|\prod_k\max\{1,|\beta_k|\}.
\]
For $|a|\leq1/2$,
\[
 \frac{|a-\beta|}{\max\{1,|\beta|\}}
 \geq \frac12\min\{1,|a-\beta|\},
\]
which proves \eqref{eq:product}. The constant-polynomial case is included by the empty-product convention.

Cauchy--Schwarz gives, on $|s|\leq1/2$,
\[
 |q'(s)|\leq
 \left(\sum_{k=1}^{\infty}k^2(1/2)^{2k-2}\right)^{1/2}
 =\sqrt{80/27}<2.
\]
Hence $|b_0|\leq |q(a)|+2R\varepsilon$ and
\[
 \norm{D_\varepsilon^{-1}b}_2
 \leq \frac{|b_0|}{\varepsilon}+1
 \leq\frac{|q(a)|}{\varepsilon}+2R+1.
\]
Splitting into $|q(a)|\leq\varepsilon$ and $|q(a)|\geq\varepsilon$ proves \eqref{eq:scaled}.
\end{proof}

\section{A lower bound at a fixed scale}\label{sec:finite}

Fix real analytic functions $f_0,w$ on $\R$, where $f_0$ has a specified simple root $\alpha_0$ and $w(\alpha_0)\ne0$. At a scale $\varepsilon>0$, set $x=\alpha_0+\varepsilon$ and consider
\begin{equation}\label{eq:family}
 f_u(t)=f_0(t)+w(t)\left[\varepsilon u_0+
                  \sum_{k=1}^{n}u_k(t-x)^k\right],
 \qquad u\in\ball_{\R^{n+1}}(0,1).
\end{equation}
Define
\[
 v_\varepsilon(t)=w(t)\bigl(\varepsilon,t-x,\dots,(t-x)^n\bigr),
 \qquad f_u(t)=f_0(t)+u\cdot v_\varepsilon(t).
\]
Assume that for all sufficiently small $\varepsilon$, every family member has a unique specified root $\alpha(u)$ in a fixed interval $J$, and that fixed positive constants $m,M,R,w_*,K$ satisfy
\begin{align}
 m\leq f'_u(t)&\leq M &&(t\in J),\label{eq:familyderiv}\\
 |\alpha(u)-x|&\leq R\varepsilon,
 & |w(\alpha(u))|&\geq w_*,\label{eq:familyroot}\\
 \norm{v'_\varepsilon(t)}_2&\leq K &&(t\in J).
 \label{eq:familyK}
\end{align}
The entire-function construction in the next section arranges these hypotheses, with $J=\R$. The argument below only needs derivative observations to be linear on the finite-dimensional family. Accordingly it also allows arbitrary real scalar linear observations at this finite-scale stage.

For the induction, define
\begin{equation}\label{eq:B}
 B_0=1,\qquad B_j=2^{j-1}\quad(j\geq1).
 \qquad\text{Then}\qquad B_j=\sum_{i=0}^{j-1}B_i\quad(j\geq1).
\end{equation}

\begin{lemma}[Lower bound at a fixed scale]\label{lem:finite}
Under \eqref{eq:familyderiv}--\eqref{eq:familyK}, there are constants $c>0$ and $\varepsilon_0>0$, independent of the algorithm, such that every adaptive deterministic algorithm with at most $n$ scalar linear observations has, for every $0<\varepsilon<\varepsilon_0$, an input from \eqref{eq:family} on which
\begin{equation}\label{eq:finitebad}
 |T_{f_u}(x)-\alpha(u)|\geq c\varepsilon^{B_n}.
\end{equation}
The constant $c$ depends only on $n,m,M,R,w_*,K$; $\varepsilon_0$ may also depend on the common scale threshold in the family hypotheses. Neither constant depends on the queries or answers in the execution, or on the particular $0<\varepsilon<\varepsilon_0$.
For derivative observations, the input can also be chosen so that no queried location equals $\alpha(u)$.
\end{lemma}

The condition on query locations will allow later corrections to vanish at these points without forcing them to vanish at the root.

\subsection{Parameters with identical observations}

We prove \cref{lem:finite} by choosing answers successively and retaining a ball of parameters consistent with them. A fixed transcript imposes affine equations on parameter space. The admissible parameters satisfying these equations are the intersection of an affine subspace with the parameter unit ball. If the first $j$ observations have been fixed, all these inputs generate the same first $j$ queries and answers: this follows inductively because query $i$ depends only on the answers before $i$.

Let $u_c$ belong to such an affine subspace, with linear direction space $E$, and consider a ball $\ball_E(u_c,r)=\{u_c+h:h\in E,\ \norm{h}_2\leq r\}$ contained in the parameter unit ball. The norm of evaluation on the remaining directions is
\begin{equation}\label{eq:sensitivity}
 a=\norm{P_Ev_\varepsilon(\alpha(u_c))}_2,
\end{equation}
where $P_E$ is orthogonal projection onto $E$. For $h\in E$ with $\norm{h}_2\leq r$, the equality
\[
 f_{u_c+h}(\alpha(u_c))=h\cdot v_\varepsilon(\alpha(u_c))
\]
and \eqref{eq:familyderiv} imply
\begin{equation}\label{eq:rootshift}
 |\alpha(u_c+h)-\alpha(u_c)|\leq \frac{a}{m}\norm{h}_2.
\end{equation}
If $a>0$, take the unit vector
$e=P_Ev_\varepsilon(\alpha(u_c))/a$. The roots at
$u_c\pm(r/2)e$ lie on opposite sides of $\alpha(u_c)$ and are each at least $ra/(2M)$ away from it. Continuity of the root along this segment shows that its image contains a real interval of length at least
\begin{equation}\label{eq:rootinterval}
 \frac{ra}{M}.
\end{equation}
We will choose the next centre on this segment. To control root motion in a smaller ball about that centre, we also need an upper bound for the projected norm there. For any $u_*$ on the segment, \eqref{eq:familyK} and \eqref{eq:rootshift} give
\begin{equation}\label{eq:sensitivitycompare}
 \norm{P_Ev_\varepsilon(\alpha(u_*))}_2
 \leq a\left(1+\frac{Kr}{2m}\right)\leq Qa,
 \qquad Q=1+\frac{K}{2m},
\end{equation}
provided $r\leq1$.

\subsection{Choice of the next answer}

Let $L_i$ be the linear observation functional chosen at step $i$ along the transcript. For a polynomial $q$ in the centred variable $s=t-x$, define
\[
 \widetilde L_i(q)=L_i\bigl(t\mapsto w(t)q(t-x)\bigr),
 \qquad
 V_{j,d}=\bigcap_{i=1}^{j}\ker(\widetilde L_i|_{\Pn_d}).
\]
Real observation functionals are extended complex-linearly on the polynomial space. Their defining matrices have real entries, so the real and complex ranks agree and these kernels have real bases.

When the $j$th query is chosen, all spaces $V_{j,d}$, $1\leq d\leq n$, are determined \emph{before} the $j$th answer is supplied. Let $S_j$ consist of the numbers $x+\zeta$ for all zeros $\zeta$ of their nonconstant Wronskians. Omit zero-dimensional spaces. For a derivative observation, also include the queried location $z_j$ in $S_j$. By \eqref{eq:Mn},
\begin{equation}\label{eq:Scard}
 |S_j|\leq H_n.
\end{equation}
We use the convention $\dist(z,\varnothing)=+\infty$. We include every degree $d\leq n$ because the degree of the eventual kernel polynomial need not be known at the current stage.

Take $\varepsilon_0\leq\min\{1,1/(2R)\}$ small enough that the family hypotheses hold. Start with centre $0$, direction space $\R^{n+1}$, and radius $r_0=1/2$. By \eqref{eq:familyroot}, the projected norm at its centre is at least $w_*\varepsilon$. Suppose after $j-1$ observations there is a ball of radius $r=r_{j-1}$ in the remaining parameter space, with projected norm at its centre satisfying
\begin{equation}\label{eq:inductivea}
 a\geq A_{j-1}\varepsilon^{B_{j-1}},
\end{equation}
and with all previously imposed distance bounds valid throughout the ball.

Choose query $j$, and form $S_j$. The interval in \eqref{eq:rootinterval} cannot be covered by the radius-$ra/[4M(H_n+1)]$ discs around its at most $H_n$ points: their intersections with the real line have total length less than half that interval's length. Hence there is a point $u_*$ on the half-radius segment such that
\begin{equation}\label{eq:centeravoid}
 \dist(\alpha(u_*),S_j)\geq\frac{ra}{4M(H_n+1)}.
\end{equation}
Supply the $j$th answer from $f_{u_*}$. If $E$ is the old direction space, the new one is
\[
 E_{\rm new}=E\cap\ker\bigl[h\mapsto L_j(h\cdot v_\varepsilon)\bigr].
\]
Set
\begin{equation}\label{eq:theta}
 \theta=\min\left\{\frac12,\frac{m}{8QM(H_n+1)}\right\},
 \qquad r_j=\theta r_{j-1}.
\end{equation}
The ball $\ball_{E_{\rm new}}(u_*,r_j)$ is contained in the previous ball, since the new centre is at most $r/2$ from the old centre and $r_j\leq r/2$. Equations \eqref{eq:rootshift} and \eqref{eq:sensitivitycompare} show that every root in this smaller ball moves from $\alpha(u_*)$ by at most
\[
 \frac{Qa}{m}\theta r\leq\frac{ra}{8M(H_n+1)}.
\]
Therefore every input in the smaller ball satisfies
\begin{equation}\label{eq:avoid}
 \dist(\alpha(u),S_j)\geq\gamma_j\varepsilon^{B_{j-1}},
 \qquad
 \gamma_j=\frac{r_{j-1}A_{j-1}}{8M(H_n+1)}>0.
\end{equation}
All earlier distance bounds are retained by containment.

\subsection{Lower bound for the projected norm}

It remains to recover a lower bound of the form \eqref{eq:inductivea} for the projected norm on the smaller ball. At an arbitrary point of that ball, choose a nonzero real polynomial $q$ of minimum degree $d$ in the common kernel of the first $j$ functionals $\widetilde L_i$. Since there are $j$ scalar constraints, $d\leq j$. Normalise its coefficient vector $b$ to norm one.

If $d=0$, the unit parameter direction corresponding to the perturbation $\varepsilon w$ lies in the common kernel and gives projected norm at least $w_*\varepsilon$, which is at least $w_*\varepsilon^{B_j}$ for $\varepsilon\leq1$.

Suppose $d\geq1$. Minimality of $d$ says that the common kernel on $\Pn_{d-1}$ is zero. Thus the observation map is injective on this $d$-dimensional space and has rank at least $d$ on $\Pn_d$. Since the nonzero polynomial $q$ belongs to its kernel on $\Pn_d$, its rank is at most $d$. The rank is therefore exactly $d$, and the kernel has dimension one. Each scalar observation increases the rank by at most one. Hence the first query indices $i_k$ at which the rank reaches $k$, for $k=1,\dots,d$, are distinct integers in $\{1,\dots,j\}$. At prefix $i_k$ the space $V_{i_k,d}$ has dimension
\[
 m_k=d+1-k.
\]
Its Wronskian is nonconstant: otherwise \eqref{eq:wrdegree} would give $V_{i_k,d}=\Pn_{m_k-1}$, which cannot contain the degree-$d$ polynomial $q$ since $m_k\leq d$.
Let $\rho_1\geq\cdots\geq\rho_d$ be the distances from $\alpha(u)-x$ to the roots of $q$, counted with multiplicity. Its $m_k$ nearest roots lie in the disc of radius $\rho_k$. Since $q\in V_{i_k,d}$, \cref{lem:wr} gives a zero of $W_{V_{i_k,d}}$ within $C(d,m_k)\rho_k$ of $\alpha(u)-x$. Its translate by $x$ belongs to $S_{i_k}$. Thus \eqref{eq:avoid} gives
\begin{equation}\label{eq:rhobound}
 \rho_k\geq\frac{\gamma_{i_k}}{C(d,d+1-k)}
                   \varepsilon^{B_{i_k-1}}.
\end{equation}
The product estimate \eqref{eq:product} yields a constant $c_j>0$ such that
\begin{equation}\label{eq:qbound}
 |q(\alpha(u)-x)|\geq c_j
       \varepsilon^{\sum_{k=1}^{d}B_{i_k-1}}
 \geq c_j\varepsilon^{B_j}.
\end{equation}
The last inequality uses $0<\varepsilon<1$ and
\[
 \{i_k-1:1\leq k\leq d\}\subseteq\{0,\dots,j-1\},
 \qquad \sum_{i=0}^{j-1}B_i=B_j.
\]
The following choice is independent of the algorithm:
\[
 C_* = \max\bigl(\{1\}\cup\{C(d,m):1\leq m\leq d\leq n\}\bigr),\qquad
 c_j=4^{-n}\prod_{i=1}^{j}\min\{1,\gamma_i/C_*\}.
\]
This choice works in \eqref{eq:qbound}: adding the unused indices only decreases a product whose factors are at most one. In particular, $0<c_j\leq1$.

To obtain an admissible parameter direction from this polynomial, we must account for the scaled constant coefficient. The vector
\[
 h_q=\frac{D_\varepsilon^{-1}b}
            {\norm{D_\varepsilon^{-1}b}_2}
\]
is a unit vector in the remaining direction space. Its function perturbation is
$w(t)q(t-x)/\norm{D_\varepsilon^{-1}b}_2$. By
\eqref{eq:scaled}, \eqref{eq:familyroot}, and \eqref{eq:qbound}, its value at the current root has magnitude at least
$A_j\varepsilon^{B_j}$, where
\begin{equation}\label{eq:Aj}
 A_0=w_*,\qquad A_j=\frac{w_*c_j}{2R+2}>0.
\end{equation}
The definitions are recursive: $A_{j-1}$ determines $\gamma_j$, which determines $c_j$ and then $A_j$. The same $A_j$ also works for $d=0$, since $A_j\leq w_*$. Taking the maximum over unit vectors in the remaining direction space gives the same lower bound for the projected norm. This holds at every point of the smaller ball and completes the induction.

\subsection{Completion of the lower bound}

After $n$ answers, the remaining ball has radius
$r_n=\tfrac12\theta^n>0$ and projected norm at its centre at least
$A_n\varepsilon^{B_n}$. Taking the two parameters displaced from its centre by $\pm(r_n/2)e$, as in \eqref{eq:rootinterval}, gives roots separated by at least
\[
 \frac{r_nA_n}{M}\varepsilon^{B_n}.
\]
Their entire transcripts, and therefore their outputs, are identical. At least one output error is at least half this root separation. Thus \eqref{eq:finitebad} holds with $c=r_nA_n/(2M)$. Every queried location was included in an avoidance set, so none equals the selected root on either input. This proves \cref{lem:finite}.

\section{A single entire counterexample}\label{sec:entire}

The function supplied by \cref{lem:finite} may depend on the scale, whereas the constant and neighbourhood in \eqref{eq:localorder} may depend on the function. To prove \cref{thm:main}, we therefore construct one entire function for which the estimate fails along a sequence of starting points. Fix $p>B_n$. As in the analytic perturbation construction of Brent, Winograd, and Wolfe \cite[Sections~III--V]{BWW}, each correction preserves all previously evaluated derivative values. Here these values come from executions at several different starting points.

\subsection{Inductive construction}

Start with $f_0(t)=t$. For each integer $s\geq1$, we inductively construct an entire function of the form
\begin{equation}\label{eq:fs}
 f_s(t)=t+e^{-t^2}P_s(t),\qquad P_s\in\R[t],
\end{equation}
and a starting point $x_s$, with output
$y_s=T_{f_s}(x_s)$. We choose each correction $f_s-f_{s-1}$ small both in real $C^1$ norm and on expanding complex discs.

The correction bounds must ensure convergence and keep every later displacement of the root small compared with the error and scale at each completed stage. Let $0<b_1\leq1/32$. After stage $s$ choose the next positive bound subject to
\begin{equation}\label{eq:budget}
 b_{s+1}\leq\min\left\{
\frac{b_s}{2},\ \frac{E_s}{32},\
 \frac{\varepsilon_s}{32}\right\},
 \qquad E_s=|y_s-\alpha_s|,
\end{equation}
where $\alpha_s$ is the unique real root of $f_s$. Positivity of $E_s$ will follow from the finite-scale lower bound. The resulting bounds satisfy $\sum_{s\geq1}b_s\leq1/16$. In particular, all intermediate functions will have derivative between $15/16$ and $17/16$ on the real line.

Suppose that the first $s-1$ functions and starting points have been chosen. Let $Z_{s-1}$ be the finite set of points queried in these executions, and let $R_{s-1}(z)$ be the largest derivative order queried at $z$ during them. Assume inductively that $f_{s-1}(z)\ne0$ for every $z\in Z_{s-1}$. The next correction must vanish at each such point to order at least $R_{s-1}(z)+1$. Define
\begin{equation}\label{eq:multiplier}
 w_s(t)=\lambda_s e^{-t^2}
           \prod_{z\in Z_{s-1}}(t-z)^{R_{s-1}(z)+1},\qquad\lambda_s>0.
\end{equation}
Since $\alpha_{s-1}$ is the unique real root, $w_s(\alpha_{s-1})\ne0$.

Choose $\lambda_s$ sufficiently small that every perturbation
\begin{equation}\label{eq:stagefamily}
 g_u(t)=w_s(t)\left[\varepsilon u_0+
                      \sum_{k=1}^{n}u_k(t-x)^k\right],
 \qquad\norm{u}_2\leq1,
\end{equation}
with $x\in[-1,1]$ and $0\leq\varepsilon\leq1$, satisfies
\begin{equation}\label{eq:smallnorms}
 \sup_{t\in\R}|g_u(t)|+\sup_{t\in\R}|g'_u(t)|
       +\sup_{|z|\leq s+2}|g_u(z)|\leq b_s.
\end{equation}
Such a choice is possible: a polynomial times $e^{-t^2}$ and its first derivative are bounded on $\R$, are bounded on each complex disc, and the indicated coefficient and starting-point sets are compact.

Set $x=\alpha_{s-1}+\varepsilon$, with $0<\varepsilon\leq1/4$. All $f_{s-1}+g_u$ are globally strictly increasing, with derivative between $1/2$ and $3/2$, and differ from $t$ by a bounded amount. Each therefore has a unique real root, denoted by $\alpha(u)$. Also $|\alpha_{s-1}|\leq1/4$, so $x\in[-1,1]$.

At the old root, Cauchy--Schwarz gives
\[
 \left|\varepsilon u_0+\sum_{k=1}^{n}u_k(-\varepsilon)^k\right|
 \leq2\varepsilon.
\]
Taking $\varepsilon=1$ and $u=(1,0,\dots,0)$ in the uniform choice \eqref{eq:smallnorms} gives $\norm{w_s}_{\infty,\R}\leq b_s$. Hence the new root moves from $\alpha_{s-1}$ by at most
$4b_s\varepsilon\leq\varepsilon/4$. Thus
\begin{equation}\label{eq:scaleroots}
 \frac34\varepsilon\leq x-\alpha(u)\leq\frac54\varepsilon.
\end{equation}
For sufficiently small $\varepsilon$, $w_s$ is bounded away from zero at all these roots. Since \eqref{eq:smallnorms} holds on the whole parameter unit ball, it also gives
\[
 \norm{v'_\varepsilon(t)}_2
 =\sup_{\norm{u}_2\leq1}|u\cdot v'_\varepsilon(t)|
 =\sup_{\norm{u}_2\leq1}|g'_u(t)|\leq b_s
 \qquad(t\in\R),
\]
which is \eqref{eq:familyK} with $K=b_s$. The multiplier $w_s$ is now fixed. Its vanishing orders at earlier query points affect the constant in the lower bound, but each new execution still imposes at most $n$ linear conditions on a polynomial of degree at most $n$. Thus \cref{lem:finite} applies with a constant $c_s>0$ and a scale threshold, both fixed before $\varepsilon_s$ is chosen.

Choose $\varepsilon_s$ below that threshold and so small that
\begin{equation}\label{eq:epschoice}
 0<\varepsilon_s<\min\{1/s,1/4\},\qquad
 \frac{c_s}{2}\varepsilon_s^{B_n}
       \geq s(2\varepsilon_s)^p.
\end{equation}
Such a choice is possible because $p>B_n$.

Apply \cref{lem:finite} to choose $u$ and set
$f_s=f_{s-1}+g_u$, $x_s=\alpha_{s-1}+\varepsilon_s$, and $y_s=T_{f_s}(x_s)$. Then
\begin{equation}\label{eq:Es}
 E_s=|y_s-\alpha_s|\geq c_s\varepsilon_s^{B_n}>0.
\end{equation}
No queried location equals $\alpha_s$. Global strict monotonicity therefore implies that $f_s$ is nonzero at every new query location. At old query locations its function values remain the previous nonzero values, because \eqref{eq:multiplier} vanishes there to the required orders. This completes the induction and permits the next choice in \eqref{eq:budget}.

\subsection{Passage to the limit}

By \eqref{eq:smallnorms}, the series of entire corrections converges uniformly on every fixed complex disc. Its limit $f$ is entire and real on $\R$. The real derivatives also converge uniformly, and the summable bounds give
$3/4\leq f'(t)\leq5/4$ on $\R$. Since $f(t)-t$ is bounded, $f$ has exactly one real zero $\alpha$.

For every fixed stage $s$,
\[
 \norm{f-f_s}_{\infty,\R}\leq\sum_{k>s}b_k\leq2b_{s+1}.
\]
The lower derivative bound $f'\geq1/2$ and \eqref{eq:budget} consequently give
\begin{equation}\label{eq:tailroot}
 |\alpha-\alpha_s|\leq4b_{s+1}
 \leq\min\{E_s/8,\varepsilon_s/8\}.
\end{equation}

Every correction after stage $s$ vanishes at every point queried in that execution to an order strictly larger than its queried derivative order. Uniform convergence on complex neighbourhoods permits passage to all derivatives. Thus the corresponding derivative values of $f$ and $f_s$ are exactly equal. Induction through the adaptive query sequence gives
\begin{equation}\label{eq:sameoutput}
 T_f(x_s)=T_{f_s}(x_s)=y_s.
\end{equation}
Equations \eqref{eq:scaleroots} and \eqref{eq:tailroot} show
\[
 \frac58\varepsilon_s\leq x_s-\alpha\leq\frac{11}{8}\varepsilon_s
 <2\varepsilon_s,
\]
so $x_s\to\alpha$ through non-root starting points. Finally,
\[
 \frac{|T_f(x_s)-\alpha|}{|x_s-\alpha|^p}
 \geq\frac{E_s/2}{(2\varepsilon_s)^p}
 \geq\frac{(c_s/2)\varepsilon_s^{B_n}}{(2\varepsilon_s)^p}
 \geq s,
\]
by \eqref{eq:Es} and \eqref{eq:epschoice}. This proves \cref{thm:main}.

\begin{remark}[Simultaneous failure of larger orders]
The real construction can use $p_s=B_n+1/s$ in place of the fixed $p$ in \eqref{eq:epschoice}. For each fixed $p>B_n$, eventually $p\geq p_s$ and $2\varepsilon_s<1$, so the resulting error ratio for exponent $p$ is still at least $s$. Thus the function and sequence in \cref{thm:main} can be chosen independently of $p>B_n$.
\end{remark}

\subsection{Evaluations in prescribed groups}\label{sec:groups}

We adapt the finite-scale argument to grouped queries, then apply the entire-function construction above.

\begin{proof}[Proof of \cref{cor:batch}]
Let $n=\sum_{j=1}^{k}\ell_j$, retain the degree-$n$ polynomial perturbation space, and define
\[
 C_0=1,\qquad C_1=\ell_1,\qquad
 C_j=(\ell_j+1)C_{j-1}\quad(2\leq j\leq k).
\]
Then $C_j=\sum_{i=1}^{j}\ell_iC_{i-1}$ for $1\leq j\leq k$.
We adapt the proof of \cref{lem:finite}, retaining one parameter ball after each group. Order the queries within each group arbitrarily. All queries in a group are determined by earlier groups' answers. Hence the union $\mathcal S_j$ of the forbidden sets for its intermediate scalar prefixes is known before any answer in group $j$ is supplied, and $|\mathcal S_j|\leq\ell_jH_n$.

To avoid this union and contract the ball once, use the constants $Q,C_*$ from the scalar proof and choose
\[
 F=nH_n,\qquad
 \theta_{\rm G}=\min\left\{\frac12,\frac{m}{8QM(F+1)}\right\},
 \qquad r_j^{\rm G}=\frac12\theta_{\rm G}^{\,j}.
\]
Starting with $A_0^{\rm G}=w_*$, define recursively
\[
 \gamma_j^{\rm G}
   =\frac{r_{j-1}^{\rm G}A_{j-1}^{\rm G}}{8M(F+1)},\qquad
 c_j^{\rm G}=4^{-n}\prod_{i=1}^{j}
       \min\{1,\gamma_i^{\rm G}/C_*\}^{\ell_i},\qquad
 A_j^{\rm G}=\frac{w_*c_j^{\rm G}}{2R+2}
 \quad(1\leq j\leq k).
\]
These constants are positive and independent of the algorithm and $\varepsilon$. Retain the scalar restriction on $\varepsilon_0$.

After group $j$, the induction retains a ball of radius $r_j^{\rm G}$ with identical observations, projected norm at least $A_j^{\rm G}\varepsilon^{C_j}$, and
\[
 \dist(\alpha(u),\mathcal S_i)
   \geq\gamma_i^{\rm G}\varepsilon^{C_{i-1}}
 \qquad(1\leq i\leq j)
\]
throughout the ball. The initial ball has radius $1/2$ and projected norm at least $w_*\varepsilon$. Before group $j$, let $r,a$ be the preceding radius and projected norm at its centre. The root-interval argument selects a centre whose root is at distance at least $ra/[4M(F+1)]$ from $\mathcal S_j$. Supply every group answer there and intersect the old direction space with all group kernels. The affine ball about this centre, with that direction space and radius $\theta_{\rm G}r$, is contained in its predecessor. The root-motion bound loses at most $ra/[8M(F+1)]$, giving the required avoidance and preserving earlier bounds.

At any point of the new ball, choose a nonzero real kernel polynomial $q$ of minimum degree, with coefficient norm one. If $q$ is nonconstant, let $h_i$ count its rank-increment indices in group $i$ in the scalar argument. Then $h_i\leq\ell_i$, and each corresponding root distance is at least $(\gamma_i^{\rm G}/C_*)\varepsilon^{C_{i-1}}$. Thus \eqref{eq:product} gives
\[
 |q(\alpha(u)-x)|
 \geq4^{-n}\prod_{i=1}^{j}
   \min\{1,\gamma_i^{\rm G}/C_*\}^{h_i}
   \varepsilon^{\sum_{i=1}^{j}h_iC_{i-1}}
 \geq c_j^{\rm G}\varepsilon^{C_j}.
\]
Adding unused factors decreases the product, and $\sum_{i=1}^{j}h_iC_{i-1}\leq C_j$. Equation \eqref{eq:scaled} gives the required projected norm, since $c_j^{\rm G}\leq1$, $C_j\geq1$, and $\varepsilon<1$. The constant-polynomial case follows from the direction $\varepsilon w$, as before.

After $k$ contractions, the identical-transcript argument gives an output error at least $c_{\rm G}\varepsilon^{C_k}$, where
\[
 c_{\rm G}=\frac{r_k^{\rm G}A_k^{\rm G}}{2M}>0.
\]
The avoidance sets include every actual query location. The entire-function construction therefore applies with $C_k=\ell_1\prod_{j=2}^{k}(\ell_j+1)$ in place of $B_n$, proving \eqref{eq:batchbound}.
\end{proof}

\paragraph{\textbf{Acknowledgments and AI declaration.}}
The formal Lean verification was carried out with assistance from ChatGPT 6. The resulting Lean proof and its problem formulation were thoroughly checked by the author and with the Lean 4 Comparator against the theorems stated in this paper, including checks that they introduce no new axioms or unproved results. I wish to thank Gonzalo Cao-Labora and George Stepaniants for their kind advice and for taking the time to explain how to set this up. ChatGPT 5.6 and 6 were also used as interactive aids in discussing and checking some of the arguments in the paper. All mathematical arguments, statements, and proofs included in the paper were reviewed and adopted by the author, who takes full responsibility for the content.

\appendix
\section{Sharpness}\label{sec:sharpness}

The upper bound in \cref{thm:main} is attained for every $n\geq2$ by the methods of Kung and Traub \cite{KT}. For completeness, we give the inverse Hermite interpolation construction to verify sharpness for \eqref{eq:localorder}. Further constructions are surveyed in \cite{PP}.

Fix $n\geq2$ and a real analytic function $f$ with a simple zero $\alpha$. Choose an interval $I$ about $\alpha$ on which $f$ has an analytic inverse $g$. Set $y_0=x$, $t_i=f(y_i)$, and evaluate $t_0$ and $f'(x)$. The linear interpolant
\[
 \mathcal H_0(t)=x+\frac{t-t_0}{f'(x)}
\]
gives $y_1=\mathcal H_0(0)$. For $1\leq j\leq n-2$, evaluate $t_j$ and let $\mathcal H_j$ be the polynomial of degree at most $j+1$ satisfying
\[
 \mathcal H_j(t_0)=x,\qquad \mathcal H_j'(t_0)=1/f'(x),\qquad
 \mathcal H_j(t_i)=y_i\quad(1\leq i\leq j).
\]
Set $y_{j+1}=\mathcal H_j(0)$. Whenever an evaluated $t_i$ is zero, return $y_i$ immediately.

We verify that these interpolants are defined for every sufficiently small nonzero $x-\alpha$ until termination. Let $e_i=y_i-\alpha$. The Hermite remainder, including $j=0$, is
\begin{equation}\label{eq:hermite}
 e_{j+1}=(-1)^{j+1}D_jt_0^2\prod_{i=1}^{j}t_i,
 \qquad D_j=g[0,t_0,t_0,t_1,\dots,t_j],
\end{equation}
where the brackets denote divided differences, with repeated nodes interpreted by derivatives. Choose a compact interval $K\subset f(I)$ with zero in its interior. The divided-difference mean-value formula \cite[Section~8, equation~(44)]{deBoor2005} gives $|D_j|\leq D$ whenever its nodes lie in $K$, where
\[
 D=\max\left\{1,\max_{0\leq j\leq n-2}
       \frac{\sup_{t\in K}|g^{(j+2)}(t)|}{(j+2)!}\right\}.
\]
This bound also holds when nodes coincide. Choose $L\geq1$ such that $|f(y)|\leq L|y-\alpha|$ near $\alpha$. Shrink $\delta>0$ so that $[\alpha-\delta,\alpha+\delta]\subset I$, its image lies in $K$, and
\[
 DL^{j+2}\delta^{j+1}<\frac12
 \qquad(0\leq j\leq n-2).
\]
For $0<|e_0|<\delta$, \eqref{eq:hermite} gives $|e_1|<|e_0|/2$. Inductively, if no root has yet been found, then for $j\geq1$,
\[
 \frac{|e_{j+1}|}{|e_j|}
 \leq DL^{j+2}|e_0|^2\prod_{i=1}^{j-1}|e_i|
 <\frac12.
\]
Every new nonroot node is therefore closer to $\alpha$ than all its predecessors. The nodes remain in $I$ and are distinct, so their function values are distinct as well. This proves local admissibility, including exact termination.

The same remainder gives $e_j=O(e_0^{2^j})$ by induction, since $e_1=O(e_0^2)$ and
\[
 2+\sum_{i=1}^{j}2^i=2^{j+1}.
\]
The output $y_{n-1}$ uses at most $n$ observations. For $f(x)=e^x-1$ at zero, $g(t)=\log(1+t)$ and the leading coefficients $e_j\sim\kappa_jx^{2^j}$ satisfy
\[
 \kappa_1=\frac12,\qquad
 \kappa_{j+1}=\frac{\prod_{i=1}^{j}\kappa_i}{j+2}.
\]
These coefficients are positive, so the order $2^{n-1}$ is attained.

\section{Complex scalar observations}\label{sec:complex}

We now allow complex starting points, complex evaluation points, and complex outputs. Each value $f^{(\nu)}(z)\in\C$, with $\nu\geq0$ an integer, counts as one observation. The inputs are holomorphic functions on open domains with a specified simple zero. As in the real case, each query must lie in the input's domain, and local order requires defined executions at every sufficiently close starting point distinct from the zero. The choice of evaluations and the processing of the data remain arbitrary and deterministic.

\begin{theorem}[Complex local order bound]\label{thm:complex}
Let $n\geq1$. For every deterministic stationary method without memory using at most $n$ complex scalar derivative observations and every real $p>B_n$, there exist an entire function $f$ and complex starting points $x_s\to\alpha$ such that
\[
 |f'(z)-1|\leq\frac14\quad(|z|\leq1),
\]
$f$ has exactly one zero $\alpha$ in the unit disc, $x_s\ne\alpha$ for every $s$, and
\[
 \frac{|T_f(x_s)-\alpha|}{|x_s-\alpha|^p}\rightarrow\infty.
\]
\end{theorem}
\begin{proof}
The polynomial Wronskian and coefficient lemmas hold over $\C$. We first obtain a disc of attainable roots in each remaining affine parameter space. The lower bound at a fixed scale then follows from the same polynomial argument. Finally, we construct a single entire function by polynomial corrections.

\paragraph{Variation of the root.}
Use the complex version of \eqref{eq:family}, with $u\in\ball_{\C^{n+1}}(0,1)$. Assume that each family member has a unique root $\alpha(u)$ in a convex disc $J$, that these roots lie in a fixed compact subdisc, and that the uniform bounds $|\alpha(u)-x|\leq R\varepsilon$ and $|w(\alpha(u))|\geq w_*>0$ from \eqref{eq:familyroot} hold. Suppose also that, for every family member,
\begin{equation}\label{eq:complexbilip}
 m|z-t|\leq|f_u(z)-f_u(t)|\leq M|z-t|\qquad(z,t\in J).
\end{equation}
Let $v=v_\varepsilon$, and assume that $\norm{v(z)-v(t)}_2\leq K|z-t|$ for $z,t\in J$. In the construction below these bounds follow by integrating $|f'_u-1|\leq1/4$ and a derivative bound along straight segments.

Consider a ball $\ball_E(u_c,r)$ in the remaining affine parameter space, contained in the parameter unit ball, with complex direction space $E$. Let $\alpha_c=\alpha(u_c)$. Define
\[
 a=\norm{P_E\overline{v(\alpha_c)}}_2.
\]
When $a>0$, set $e=P_E\overline{v(\alpha_c)}/a$. The projection uses the Hermitian inner product, taken linear in its first argument, whereas the perturbation remains the complex bilinear expression $u\cdot v$. Thus $e\cdot v(\alpha_c)=a$. For $h\in E$ with $\norm{h}_2\leq r$, the root-shift estimate
\[
 |\alpha(u_c+h)-\alpha_c|\leq a\norm{h}_2/m
\]
follows from \eqref{eq:complexbilip} exactly as in the real case.

Choose the initial radius $r_0=\min\{1/2,2M/K\}$, with $2M/K=+\infty$ when $K=0$, and suppose $0<r\leq r_0$. As the parameter varies over $u_c+\{\xi e:|\xi|\leq r/2\}$, its root ranges over a set containing the disc
\begin{equation}\label{eq:rootimagedisk}
 \ball\left(\alpha_c,\frac{ra}{4M}\right).
\end{equation}
Indeed, for any target $t$ in this disc, set
\[
 \xi=-\frac{f_{u_c}(t)}{e\cdot v(t)}.
\]
The denominator has magnitude at least
$a-K|t-\alpha_c|\geq a/2$, and the numerator has magnitude at most $M|t-\alpha_c|$. Consequently $|\xi|\leq r/2$ and $f_{u_c+\xi e}(t)=0$. Thus $t$ is that input's selected root. For sufficiently small $\varepsilon$ all discs in \eqref{eq:rootimagedisk} lie in $J$: uniformly, $a\leq\norm{v(\alpha_c)}_2=O(\varepsilon)$, and the roots have a fixed positive distance from the boundary of $J$.

\paragraph{Avoidance and projected norm.}
Form the same sets $S_j$ before each answer, including the current query location when present. The disc \eqref{eq:rootimagedisk} cannot be covered by its at most $H_n$ forbidden discs of radius
\[
 \delta=\frac{ra}{8M(H_n+1)};
\]
the sum of their areas is smaller than the area of the disc in \eqref{eq:rootimagedisk}. Choose $u_*$ within distance $r/2$ of the previous centre whose root avoids $S_j$ by $\delta$, and supply its answer. With $Q=1+K/(2m)$, choose
\[
 \theta=\min\left\{\frac12,\frac{m}{16QM(H_n+1)}\right\},
 \qquad r_j=\theta r_{j-1}.
\]
The ball of radius $\theta r$ in the new affine parameter space has all its roots at least $\delta/2$ from $S_j$. The polynomial argument in \cref{sec:finite} applies over $\C$, with
\[
 \gamma_j=\frac{r_{j-1}A_{j-1}}{16M(H_n+1)}
\]
and the same recursive definitions of $c_j,A_j$. It yields a final projected norm at least $A_n\varepsilon^{B_n}$. Opposite points of the final disc of attainable roots correspond to two inputs with identical observations and roots separated by at least $r_nA_n\varepsilon^{B_n}/(2M)$. One error is therefore at least $c\varepsilon^{B_n}$, where $c>0$ is independent of the algorithm and $\varepsilon$.

\paragraph{An entire limit.}
Only bounds on complex discs are needed for the limiting construction, so polynomial corrections suffice. Start with $f_0(z)=z$ and keep each $f_s$ polynomial, with $\alpha_s$ denoting its unique zero in the unit disc. For the past complex query nodes, use
\[
 w_s(z)=\lambda_s\prod_{\zeta\in Z_{s-1}}
                      (z-\zeta)^{R_{s-1}(\zeta)+1}
\]
in place of \eqref{eq:multiplier}. Choose $\lambda_s>0$ so that all perturbations of the form \eqref{eq:stagefamily}, now with complex $u,x$ satisfying $\norm{u}_2\leq1$, $|x|\leq1$, and real $0\leq\varepsilon\leq1$, obey
\[
 \sup_{|z|\leq s+2}|g_u(z)|+
 \sup_{|z|\leq1}|g'_u(z)|\leq b_s.
\]
Choose the bounds $b_s$ according to \eqref{eq:budget}. At stage $s$, let $f_u=f_{s-1}+g_u$. Throughout the construction,
$|f_u(z)-z|\leq1/16$ and $|f'_u(z)-1|\leq1/16$ on the unit disc. Rouch\'e's theorem gives exactly one zero there, denoted by $\alpha(u)$, and the equation $f_u(\alpha(u))=0$ shows $|\alpha(u)|\leq1/16$. In particular, the roots remain in a fixed compact subdisc. For $z,t$ in the unit disc,
\[
 f_u(z)-f_u(t)=(z-t)+\int_{[t,z]}(f'_u(\zeta)-1)\,\mathrm d\zeta.
\]
This gives \eqref{eq:complexbilip} with $m=1/2$ and $M=3/2$.

Every old query node inside the unit disc has a nonzero function value that is preserved by all subsequent corrections: when first queried it was not the unique selected root. An old node outside that disc cannot equal any selected root. Hence $w_s(\alpha_{s-1})\ne0$.

Set $x=\alpha_{s-1}+\varepsilon$. At the old root the perturbation is at most $2b_s\varepsilon$, so \eqref{eq:complexbilip} yields a root shift at most $4b_s\varepsilon\leq\varepsilon/4$. Consequently
\[
 \frac34\varepsilon\leq|x-\alpha(u)|\leq\frac54\varepsilon.
\]
With $w_s$ fixed, the complex finite-scale argument supplies a constant $c_s>0$ and a scale threshold. Choose $\varepsilon_s$ below this threshold and satisfying \eqref{eq:epschoice}, and then select $f_s$ and $x_s$ as in \cref{sec:entire}. Uniform convergence on every complex disc gives an entire limit, preserves every observed derivative value, and yields
\[
 |\alpha-\alpha_s|\leq4b_{s+1}
 \leq\min\{E_s/8,\varepsilon_s/8\}.
\]
The same induction through the adaptive queries preserves the outputs. The final error-ratio estimate follows as in the real proof, using absolute values in \eqref{eq:scaleroots}. This proves the theorem.
\end{proof}

\begingroup
\small
\bibliographystyle{amsplain}
\bibliography{kung_traub_solution}
\endgroup
\end{document}